\RequirePackage{ifpdf}
\ifpdf % We are running pdfTeX in pdf mode
\documentclass[pdftex]{sigma}
\else
\documentclass{sigma}
\fi

\begin{document}

\renewcommand{\thefootnote}{$\star$}

\renewcommand{\PaperNumber}{082}

\FirstPageHeading

\ShortArticleName{Some Orthogonal Polynomials in Four Variables}

\ArticleName{Some Orthogonal Polynomials in Four Variables\footnote{This paper is a contribution to the Special
Issue on Dunkl Operators and Related Topics. The full collection
is available at
\href{http://www.emis.de/journals/SIGMA/Dunkl_operators.html}{http://www.emis.de/journals/SIGMA/Dunkl\_{}operators.html}}}

\Author{Charles F. DUNKL}

\AuthorNameForHeading{C.F. Dunkl}

\Address{Department of Mathematics, University of Virginia, Charlottesville, VA 22904-4137, USA}

\Email{\href{mailto:cfd5z@virginia.edu}{cfd5z@virginia.edu}}
\URLaddress{\url{http://people.virginia.edu/~cfd5z/}}

\ArticleDates{Received October 14, 2008, in f\/inal form November 24,
2008; Published online November 29, 2008}

\Abstract{The symmetric group on 4 letters has the ref\/lection group $D_{3}$ as an
isomorphic image. This fact follows from the coincidence of the root systems
$A_{3}$ and $D_{3}$. The isomorphism is used to construct an orthogonal basis
of polynomials of 4 variables with~2 parameters. There is an associated
quantum Calogero--Sutherland model of 4 identical particles on the line.}


\Classification{33C52; 05E35; 37J35}

\section{Introduction}

The symmetric group on $N$ letters acts naturally on $\mathbb{R}^{N}$ (for $N=2,3,\ldots$) but not irreducibly, because the vector $\left(
1,1,\ldots,1\right)  $ is f\/ixed. However the important basis consisting of
nonsymmetric Jack polynomials is def\/ined for $N$ variables and does not behave
well under restriction to the orthogonal complement of $\left(  1,1,\ldots
,1\right)  $, in general. In this paper we consider the one exception to this
situation, occurring when $N=4$. In this case there is a coordinate system,
essentially the $4\times4$ Hadamard matrix, which allows a dif\/ferent basis of
polynomials, derived from the type-$B$ nonsymmetric Jack polynomials for the
subgroup $D_{3}$ of the octahedral group $B_{3}$. We will construct an
orthogonal basis for the $L^{2}$-space of the measure%
\[
\prod_{1\leq i<j\leq4}\left\vert x_{i}-x_{j}\right\vert ^{2\kappa}\left\vert
x_{1}+x_{2}+x_{3}+x_{4}\right\vert ^{2\kappa^{\prime}}\exp\left(  -\frac{1}%
{2}\sum_{i=1}^{4}x_{i}^{2}\right)  dx
\]
on $\mathbb{R}^{4}$, with $\kappa,\kappa^{\prime}>0$.

We will use the following notations: $\mathbb{N}_{0}$ denotes the set of nonnegative integers; $\mathbb{N}_{0}^{N}$ is the set of compositions (or multi-indices), if $\alpha=\left(
\alpha_{1},\ldots,\alpha_{N}\right)  \in \mathbb{N}_{0}^{N}$ then $\left\vert \alpha\right\vert :=\sum_{i=1}^{N}\alpha_{i}$ and
the length of $\alpha$ is $\ell\left(  \alpha\right)  :=\max\left\{
i:\alpha_{i}>0\right\}  $. Let $\mathbb{N}_{0}^{N,+}$ denote the subset of partitions, that is, $\lambda\in
\mathbb{N}
_{0}^{N}$ and $\lambda_{i}\geq\lambda_{i+1}$ for $1\leq i<N$. For $\alpha\in \mathbb{N}_{0}^{N}$ and $x\in
\mathbb{R}^{N}$ let $x^{\alpha}=\prod_{i=1}^{N}x_{i}^{\alpha_{i}}$, a monomial of degree~$\left\vert \alpha\right\vert $; the space of polynomials is $\mathcal{P}%
=\mathrm{span}_{\mathbb{R}}\left\{  x^{\alpha}:\alpha\in
\mathbb{N}_{0}^{N}\right\}  $. For $x,y\in
\mathbb{R}^{N}$ the inner product is $\left\langle x,y\right\rangle :=\sum_{i=1}^{N}x_{i}y_{i}$, and $\left\vert x\right\vert :=\left\langle x,x\right\rangle
^{1/2}$; also $x^{\bot}:=\left\{  y:\left\langle x,y\right\rangle =0\right\}
$. The cardinality of a set $E$ is denoted by $\#E$.

Consider the elements of $S_{N}$ as permutations on $\{1,2,\ldots,N\}$. For
$x\in\mathbb{R}^{N}$ and $w\in\mathcal{S}_{N}$ let $\left(  xw\right)
_{i}:=x_{w\left(  i\right)  }$ for $1\leq i\leq N$ and extend this action to
polynomials by $\left(  wf\right)  \left(  x\right)  =f\left(  xw\right)  $.
Monomials transform to monomials: $w\left(  x^{\alpha}\right)  :=x^{w\alpha}$
where $\left(  w\alpha\right)  _{i}:=\alpha_{w^{-1}\left(  i\right)  }$ for
$\alpha\in\mathbb{N}_{0}^{N}$. (Consider $x$ as a row vector, $\alpha$ as a
column vector, and $w$ as a permutation matrix, with $1$'s at the $\left(
w\left(  j\right)  ,j\right)  $ entries.) For $1\leq i\leq N$ and
$f\in\mathcal{P}$ the Dunkl operators are%
\[
\mathcal{D}_{i}f\left(  x\right)  :=\frac{\partial}{\partial x_{i}}f\left(
x\right)  +\kappa\sum_{j\neq i}\frac{f\left(  x\right)  -f\left(  x\left(
i,j\right)  \right)  }{x_{i}-x_{j}},
\]
and
\[
\mathcal{U}_{i}f\left(  x\right)  :=\mathcal{D}_{i}\left(  x_{i}f\left(
x\right)  \right)  -\kappa\sum_{j=1}^{i-1}\left(  j,i\right)  f\left(
x\right)  .
\]

Then $\mathcal{U}_{i}\mathcal{U}_{j}=\mathcal{U}_{j}\mathcal{U}_{i}$ for
$1\leq i,j\leq N$ and these operators are self-adjoint for the following
pairing
\[
\left\langle f,g\right\rangle _{\kappa}:=f\left(  \mathcal{D}_{1}%
,\ldots,\mathcal{D}_{N}\right)  g\left(  x\right)  |_{x=0}.
\]
This satisf\/ies $\left\langle f,g\right\rangle _{\kappa}=\left\langle
g,f\right\rangle _{\kappa}=\left\langle wf,wg\right\rangle _{\kappa}$ for
$f,g\in\mathcal{P}$ and $w\in\mathcal{S}_{N}$; furthermore $\left\langle
f,f\right\rangle _{\kappa}>0$ when $f\neq0$ and $\kappa\geq0$. The operators
$\mathcal{U}_{i}$ have the very useful property of acting as triangular
matrices on the monomial basis furnished with a certain partial order. However
the good properties depend completely on the use of $%
%TCIMACRO{\U{211d} }%
%BeginExpansion
\mathbb{R}
%EndExpansion
^{N}$ even though the group $\mathcal{S}_{N}$ acts irreducibly on $\left(
1,1,\ldots,1\right)  ^{\bot}$. We suggest that an underlying necessity for the
existence of an analog of $\left\{  \mathcal{U}_{i}\right\}  $ for any
ref\/lection group $W$ is the existence of a $W$-orbit in which any two points
are orthogonal or antipodal (as in the analysis of the hyperoctahedral group
$B_{N}$). This generally does not hold for the action of $\mathcal{S}_{N}$ on
$\left(  1,\ldots,1\right)  ^{\bot}$. We consider the exceptional case $N=4$
and exploit the isomorphism between $\mathcal{S}_{4}$ and the group of type
$D_{3}$, that is, the subgroup of $B_{3}$ whose simple roots are $\left(
1,-1,0\right)$, $\left(  0,1,-1\right)$, $\left(  0,1,1\right)  $. We map these
root vectors to the simple roots $\left(  0,1,-1,0\right)$, $\left(
0,0,1,-1\right)$, $\left(  1,-1,0,0\right)  $ of $\mathcal{S}_{4}$, in the same
order. This leads to the linear isometry%
\begin{gather}
y_{1}    =\tfrac{1}{2}\left(  x_{1}+x_{2}-x_{3}-x_{4}\right), \nonumber\\
y_{2}    =\tfrac{1}{2}\left(  x_{1}-x_{2}+x_{3}-x_{4}\right), \nonumber\\
y_{3}    =\tfrac{1}{2}\left(  x_{1}-x_{2}-x_{3}+x_{4}\right), \nonumber\\
y_{0}    =\tfrac{1}{2}\left(  x_{1}+x_{2}+x_{3}+x_{4}\right)  .\label{y2x}
\end{gather}
Consider the group $D_{3}$ acting on $\left(  y_{1},y_{2},y_{3}\right)  $ and
use the type-$B_{3}$ Dunkl operators with the parameter $\kappa^{\prime}=0$
(associated with the class of sign-changes $y_{i}\mapsto-y_{i}$ which are not
in $D_{3}$). Let~$\sigma_{ij}$, $\tau_{ij}$ denote the ref\/lections in
$y_{i}-y_{j}=0$, $y_{i}+y_{j}=0$ respectively. Then for $i=1,2,3$ let%
\begin{gather*}
\mathcal{D}_{i}^{B}f\left(  y\right)      =\frac{\partial}{\partial y_{i}%
}f\left(  y\right)  +\kappa\sum_{j=1,j\neq i}^{3}\left(  \frac{f\left(
y\right)  -f\left(  y\sigma_{ij}\right)  }{y_{i}-y_{j}}+\frac{f\left(
y\right)  -f\left(  y\tau_{ij}\right)  }{y_{i}+y_{j}}\right)  ,\\
\mathcal{U}_{i}^{B}f\left(  y\right)      =\mathcal{D}_{i}^{B}\left(
y_{i}f\left(  y\right)  \right)  -\kappa\sum_{1\leq j<i}\left(  \sigma
_{ij}+\tau_{ij}\right)  f\left(  y\right)  .
\end{gather*}
The operators $\left\{  \mathcal{U}_{i}^{B}\right\}  $ commute pairwise and
are self-adjoint for the usual inner product. The simultaneous eigenvectors
are expressed in terms of nonsymmetric Jack polynomials with argument $\left(
y_{1}^{2},y_{2}^{2},y_{3}^{2}\right)  $. In the sequel we consider polynomials
with arguments $x$ or $y$ with the convention that $y$ is given in terms of
$x$ by equation~(\ref{y2x}).

\section{Nonsymmetric Jack polynomials}

\textit{Nonsymmetric Jack polynomials} (NSJP) are the simultaneous
eigenfunctions of $\left\{  \mathcal{U}_{i}\right\}  _{i=1}^{N}$. We consider
the formulae for arbitrary $N$ since there is really no simplif\/ication for
$N=3$.

\begin{definition}
For $\alpha\in\mathbb{N}_{0}^{N}$, let $\alpha^{+}$ denote the unique
partition such that $\alpha^{+}=w\alpha$ for some $w\in S_{N}$. For
$\alpha,\beta\in\mathbb{N}_{0}^{N}$ the partial order $\alpha\succ\beta$
($\alpha$ dominates $\beta$) means that $\alpha\neq\beta$ and $\sum
_{i=1}^{j}\alpha_{i}\geq\sum_{i=1}^{j}\beta_{i}$ for $1\leq
j\leq N$; $\alpha\vartriangleright\beta$ means that $\left\vert \alpha
\right\vert =\left\vert \beta\right\vert $ and either $\alpha^{+}\succ
\beta^{+}$ or $\alpha^{+}=\beta^{+}$ and $\alpha\succ\beta$.
\end{definition}

For example $\left(  2,6,4\right)  \vartriangleright\left(  5,4,3\right)
\vartriangleright\left(  3,4,5\right)  $. When acting on the monomial basis
 $\big\{  x^{\alpha}:\alpha\in
\mathbb{N}
_{0}^{N}$, $\left\vert \alpha\right\vert =n\big\}  $ for $n\in\mathbb{N}_{0}$
the operators $\mathcal{U}_{i}$ have on-diagonal coef\/f\/icients given by the
following functions on $\mathbb{N}_{0}^{N}$.

\begin{definition}
For $\alpha\in\mathbb{N}_{0}^{N}$ and $1\leq i\leq N$ let
\begin{gather*}
r\left(  \alpha,i\right)      :=\#\left\{  j:\alpha_{j}>\alpha_{i}\right\}
+\#\left\{  j:1\leq j\leq i,\alpha_{j}=\alpha_{i}\right\}  ,\\
\xi_{i}\left(  \alpha\right)      :=\left(  N-r\left(  \alpha,i\right)
\right)  \kappa+\alpha_{i}+1.
\end{gather*}
\end{definition}

Clearly for a f\/ixed $\alpha\in\mathbb{N}_{0}^{N}$ the values $\left\{
r\left(  \alpha,i\right)  :1\leq i\leq N\right\}  $ consist of all of
$\left\{  1,\ldots,N\right\}  $; let $w$ be the inverse function of $i\mapsto
r\left(  \alpha,i\right)  $ so that $w\in\mathcal{S}_{N}$, $r\left(
\alpha,w\left(  i\right)  \right)  =i$ and $\alpha^{+}=w\alpha\ $(note that
$\alpha\in\mathbb{N}_{0}^{N,+}$ if and only if $r\left(  \alpha,i\right)  =i$
for all $i$). Then%
\[
\mathcal{U}_{i}x^{\alpha}=\xi_{i}\left(  \alpha\right)  x^{\alpha}%
+q_{\alpha,i}\left(  x\right)
\]
where $q_{\alpha,i}\left(  x\right)  $ is a sum of terms $\pm\kappa x^{\beta}$
with $\alpha\vartriangleright\beta$.

\begin{definition}
For $\alpha\in\mathbb{N}_{0}^{N}$, let$\ \zeta_{\alpha}$ denote the $x$-monic
simultaneous eigenfunction (NSJP), that is, $\mathcal{U}_{i}\zeta_{\alpha}%
=\xi_{i}\left(  \alpha\right)  \zeta_{\alpha}$ for $1\leq i\leq N$ and
\[
\zeta_{\alpha}=x^{\alpha}+\sum\limits_{\alpha\vartriangleright\beta}%
A_{\beta\alpha}x^{\beta},
\]
with coef\/f\/icients $A_{\beta\alpha}\in\mathbb{Q}\left(  \kappa\right)  $,
rational functions of $\kappa$.
\end{definition}

There are norm formulae for the pairing $\left\langle \cdot,\cdot\right\rangle
_{\kappa}$.
Suppose $\alpha\in\mathbb{N}_{0}^{N}$ and $\ell\left(  \alpha\right)  =m$; the
\textit{Ferrers diagram} of $\alpha$ is the set $\left\{  \left(
i,j\right)  :1\leq i\leq m,0\leq j\leq\alpha_{i}\right\}  .$ For each node
$\left(  i,j\right)  $ with $1\leq j\leq\alpha_{i}$   there are two special
subsets of the Ferrers diagram, the \textit{arm} $\left\{  \left(  i,l\right)
:j<l\leq\alpha_{i}\right\}  $ and the \textit{leg} $\left\{  \left(
l,j\right)  :l>i,j\leq\alpha_{l}\leq\alpha_{i}\right\}  \cup\left\{  \left(
l,j-1\right)  :l<i,j-1\leq\alpha_{l}<\alpha_{i}\right\}  $. The node itself,
the arm and the leg make up the \textit{hook}. (For the case of partitions the
nodes $\left(  i,0\right)  $ are customarily omitted from the Ferrers
diagram.) The cardinality of the leg is called the leg-length, formalized by
the following:

\begin{definition}
For $\alpha\in\mathbb{N}_{0}^{N}$, $1\leq i\leq\ell\left(  \alpha\right)  $ and
$1\leq j\leq\alpha_{i}$ the leg-length is%
\[
L\left(  \alpha;i,j\right)  :=\#\left\{  l:l>i,j\leq\alpha_{l}\leq\alpha
_{i}\right\}  +\#\left\{  l:l<i,j\leq\alpha_{l}+1\leq\alpha_{i}\right\}  .
\]
For $t\in\mathbb{Q}\left(  \kappa\right)  $ the \textit{hook-length }and the
hook-length product for $\alpha$ are given by
\begin{gather*}
h\left(  \alpha,t;i,j\right)      :=\left(  \alpha_{i}-j+t+\kappa L\left(
\alpha;i,j\right)  \right), \\
h\left(  \alpha,t\right)      :=\prod_{i=1}^{\ell\left(  \alpha\right)  }%
\prod_{j=1}^{\alpha_{i}}h\left(  \alpha,t;i,j\right)  ,
\end{gather*}
and for $\lambda\in\mathbb{N}_{0}^{N,+}$ and $t\in\mathbb{Q}\left(
\kappa\right)  $ the generalized Pochhammer symbol is
\[
\left(  t\right)  _{\lambda}:=\prod\limits_{i=1}^{N}\prod\limits_{j=0}
^{\lambda_{i}-1}\left(  t-\left(  i-1\right)  \kappa+j\right)  .
\]
\end{definition}

(The product over $j$ is an ordinary Pochhammer symbol.)

\begin{proposition}
For $\alpha,\beta\in
\mathbb{N}
_{0}^{N}$, the following orthogonality and norm formula holds:%
\[
\left\langle \zeta_{\alpha},\zeta_{\beta}\right\rangle _{\kappa}%
=\delta_{\alpha\beta}\left(  N\kappa+1\right)  _{\alpha^{+}}\frac{h\left(
\alpha,1\right)  }{h\left(  \alpha,\kappa+1\right)  }.
\]
\end{proposition}

Details can be found in the book by Xu and the author \cite[Chapter~8]{DX}, the
concept of leg-length and its use in the norm formula is due to Knop and Sahi~\cite{KS}. The (symmetric) Jack polynomial with leading term $x^{\lambda}$ for
$\lambda\in \mathbb{N}_{0}^{N,+}$ is obtained by symmetrizing $\zeta_{\lambda}$. The coef\/f\/icients
involve, for $\alpha\in \mathbb{N}
_{0}^{N}$, $\varepsilon=\pm1$:%
\[
\mathcal{E}_{\varepsilon}\left(  \alpha\right)  :=\prod_{i<j,\alpha_{i}%
<\alpha_{j}}\left(  1+\frac{\varepsilon\kappa}{\left(  r\left(  \alpha
,i\right)  -r\left(  \alpha,j\right)  \right)  \kappa+\alpha_{j}-\alpha_{i}%
}\right)  ,
\]
in fact, \cite[Lemma~3.10]{D1},
\begin{gather*}
h\left(  \alpha,\kappa+1\right)      =\mathcal{E}_{1}\left(  \alpha\right)
h\left(  \alpha^{+},\kappa+1\right)  ,\\
h\left(  \alpha^{+},1\right)      =h\left(  \alpha,1\right)  \mathcal{E}%
_{-1}\left(  \alpha\right)  ,
\end{gather*}
for $\alpha\in
\mathbb{N}
_{0}^{N}$. Then
\begin{gather*}
j_{\lambda}     =\sum_{\alpha^{+}=\lambda}\mathcal{E}_{-1}\left(
\alpha\right)  \zeta_{\alpha},\\
\left\langle j_{\lambda},j_{\lambda}\right\rangle _{\kappa}     =\#\left\{
\alpha:\alpha^{+}=\lambda\right\}  \frac{\left(  N\kappa+1\right)  _{\lambda
}h\left(  \lambda,1\right)  }{\mathcal{E}_{1}\left(  \lambda^{R}\right)
h\left(  \lambda,\kappa+1\right)  },
\end{gather*}
where $\lambda_{i}^{R}=\lambda_{N+1-i}$ for $1\leq i\leq N$ (the reverse of
$\lambda$). Note $\left\{  \alpha:\alpha^{+}=\lambda\right\}  =\left\{
w\lambda:w\in\mathcal{S}_{N}\right\}  $.

\section[The groups $\mathcal{S}_{4}$ and $D_{3}$]{The groups $\boldsymbol{\mathcal{S}_{4}}$ and $\boldsymbol{D_{3}}$}

By using the $x\leftrightarrow y$ correspondence (equation (\ref{y2x})) we
obtain operators which behave well on $\left(  1,\ldots,1\right)  ^{\bot}$.
Here are the lists of ref\/lections in corresponding order:
\begin{gather*}
   \left[  \sigma_{12},\tau_{12},\sigma_{13},\tau_{13},\sigma_{23},\tau
_{23}\right]  ,\\
   \left[  \left(  23\right)  ,\left(  14\right)  ,\left(  24\right)  ,\left(
13\right)  ,\left(  34\right)  ,\left(  12\right)  \right]  .
\end{gather*}
The following orthonormal basis is used in the directional derivatives:
\begin{gather*}
v_{0}    =\tfrac{1}{2}\left(  1,1,1,1\right)  ,\\
v_{1}    =\tfrac{1}{2}\left(  1,1,-1,-1\right)  ,\\
v_{2}    =\tfrac{1}{2}\left(  1,-1,1,-1\right)  ,\\
v_{3}    =\tfrac{1}{2}\left(  1,-1,-1,1\right)  .
\end{gather*}
That is, $y_{i}=\left\langle x,v_{i}\right\rangle $ and $\frac{\partial
}{\partial y_{i}}=\sum_{j=1}^{4}\left(  v_{i}\right)  _{j}\frac{\partial
}{\partial x_{j}}$ for $0\leq i\leq3$. Note that $\left\{  \pm v_{1},\pm
v_{2},\pm v_{3}\right\}  $ is an octahedron and an $\mathcal{S}_{4}$-orbit.
Then
\begin{gather*}
\mathcal{D}_{1}^{B}f\left(  x\right)      =\sum_{j=1}^{4}\left(  v_{1}\right)
_{j}\frac{\partial f\left(  x\right)  }{\partial x_{j}}
  +\kappa\left(  \frac{1-\left(  23\right)  }{x_{2}-x_{3}}+\frac{1-\left(
14\right)  }{x_{1}-x_{4}}+\frac{1-\left(  24\right)  }{x_{2}-x_{4}}%
+\frac{1-\left(  13\right)  }{x_{1}-x_{3}}\right)  f\left(  x\right)  ,
\end{gather*}
and similar expressions hold for $\mathcal{D}_{2}^{B}$, $\mathcal{D}_{3}^{B}$.
Furthermore%
\begin{gather*}
\mathcal{U}_{1}^{B}f\left(  x\right)      =\mathcal{D}_{1}^{B}\left(
\left\langle v_{1},x\right\rangle f\left(  x\right)  \right)  ,\\
\mathcal{U}_{2}^{B}f\left(  x\right)      =\mathcal{D}_{2}^{B}\left(
\left\langle v_{2},x\right\rangle f\left(  x\right)  \right)  -\kappa\left(
\left(  14\right)  +\left(  23\right)  \right)  f\left(  x\right)  ,\\
\mathcal{U}_{3}^{B}f\left(  x\right)      =\mathcal{D}_{3}^{B}\left(
\left\langle v_{3},x\right\rangle f\left(  x\right)  \right)  -\kappa\left(
\left(  12\right)  +\left(  13\right)  +\left(  24\right)  +\left(  34\right)
\right)  f\left(  x\right)  .
\end{gather*}
For a subset $E\subset\left\{  1,2,3\right\}  $ let $y_{E}=\prod_{i\in E}%
y_{i}$, also let $E_{0}=\varnothing$ and $E_{k}=\left\{  1,\ldots,k\right\}  $
for $k=1,2,3$. The simultaneous eigenfunctions are of the form $y_{E}f\left(
y^{2}\right)  $ where $y^{2}:=\left(  y_{1}^{2},y_{2}^{2},y_{3}^{2}\right)  $
and when $E=E_{k}$ with $0\leq k\leq3$ they are directly expressed as NSJP's
(for $%
%TCIMACRO{\U{211d} }%
%BeginExpansion
\mathbb{R}
%EndExpansion
^{3}$). The following is the specialization to $\kappa^{\prime}=0$ of the
type-$B$ result from \cite[Corollary~9.3.3, p.~342]{DX}.

\begin{proposition}
Suppose $\alpha\in
\mathbb{N}
_{0}^{3}$ and $k=0,1,2,3$, then for $1\leq i\leq k$%
\[
\mathcal{U}_{i}^{B}y_{E_{k}}\zeta_{\alpha}\left(  y^{2}\right)  =2\xi
_{i}\left(  \alpha\right)  y_{E_{k}}\zeta_{\alpha}\left(  y^{2}\right)  ,
\]
and for $k<i\leq3$%
\[
\mathcal{U}_{i}^{B}y_{E_{k}}\zeta_{\alpha}\left(  y^{2}\right)  =\left(
2\xi_{i}\left(  \alpha\right)  -1\right)  y_{E_{k}}\zeta_{\alpha}\left(
y^{2}\right)  .
\]
\end{proposition}

The polynomial $y_{E_{k}}\zeta_{\alpha}\left(  y^{2}\right)  $ is labeled by
$\beta\in%
%TCIMACRO{\U{2115} }%
%BeginExpansion
\mathbb{N}
%EndExpansion
_{0}^{3}$ where $\beta_{i}=2\alpha_{i}+1$ for $1\leq i\leq k$ and $\beta
_{i}=2\alpha_{i}$ for $k<i\leq3$. The dif\/ference $\beta-\alpha\in%
%TCIMACRO{\U{2115} }%
%BeginExpansion
\mathbb{N}
%EndExpansion
_{0}^{3}$ and appears in the norm formula (the result for the pairing $\left(
f,g\right)  \mapsto f\left(  \mathcal{D}_{1}^{B},\mathcal{D}_{2}%
^{B},\mathcal{D}_{3}^{B}\right)  g\left(  y\right)  |_{y=0}$ applies because
of the isomorphism).

\begin{proposition}
Suppose $\beta\in%
%TCIMACRO{\U{2115} }%
%BeginExpansion
\mathbb{N}
%EndExpansion
_{0}^{3}$ and $\beta_{i}$ is odd for $1\leq i\leq$ $k$ and is even otherwise,
then for $\alpha_{i}=\left\lfloor \frac{\beta_{i}}{2}\right\rfloor $, $1\leq
i\leq3$%
\[
\left\langle y_{E_{k}}\zeta_{\alpha}\left(  y^{2}\right)  ,y_{E_{k}}%
\zeta_{\alpha}\left(  y^{2}\right)  \right\rangle _{\kappa}=2^{\left\vert
\beta\right\vert }\left(  3\kappa+1\right)  _{\alpha^{+}}\left(  2\kappa
+\frac{1}{2}\right)  _{\left(  \beta-\alpha\right)  ^{+}}\frac{h\left(
\alpha,1\right)  }{h\left(  \alpha,\kappa+1\right)  }.
\]
\end{proposition}

(The formulae in \cite[Chapter~9]{DX} are given for the $p$-monic polynomials,
here we use the $x$-monic type, see \cite[pp.~323--324]{DX}). There is an
evaluation formula for $\zeta_{\alpha}\left(  1,1,1\right)  $ which provides
the value at $x=\left(  2,0,0,0\right)  $, corresponding to $y=\left(
1,1,1,1\right)  $. Indeed for $\alpha\in \mathbb{N}_{0}^{3}$ (see \cite[p.~324]{DX})
\[
\zeta_{\alpha}\left(  1,1,1\right)  =\frac{\left(  3\kappa+1\right)
_{\alpha^{+}}}{h\left(  \alpha,\kappa+1\right)  }.
\]

For any point $\left(  \pm2,0,0,0\right)  w$ with $w\in\mathcal{S}_{4}$ the
corresponding $y$ satisf\/ies $y_{i}=\pm1$ for $1\leq i\leq3$, so that
$y^{2}=\left(  1,1,1\right)  $. For any other subset $E\subset\left\{
1,2,3\right\}  $ with $\#E=k$ let $w\in\mathcal{S}_{3}$ be such that $w\left(
i\right)  \in E$ for $1\leq i\leq k$, $1\leq i<j\leq k$ or $k<i<j\leq3$
implies $w\left(  i\right)  <w\left(  j\right)  $ (that is, $w$ preserves
order on $\left\{  1,\ldots,k\right\}  $ and on $\left\{  k+1,\ldots
,3\right\}  $). Here is the list of sets with corresponding permutations
$\left(  w\left(  i\right)  \right)  _{i=1}^{3}$:%
\begin{alignat*}{3}
& E    =\left\{  2\right\}  ,\qquad && w=\left(  2,1,3\right)  ,& \\
& E    =\left\{  3\right\}  ,\qquad && w=\left(  3,1,2\right)  ,& \\
& E    =\left\{  1,3\right\}  ,\qquad && w=\left(  1,3,2\right)  ,& \\
& E    =\left\{  2,3\right\}  ,\qquad && w=\left(  2,3,1\right)  .&
\end{alignat*}
Then (letting $w$ act on $y$) $wy_{E_{k}}=y_{E}$ and for $\alpha\in \mathbb{N}_{0}^{3}$ the polynomial $w\left(  y_{E_{k}}\zeta_{\alpha}\left(
y^{2}\right)  \right)  $ is a~simultaneous eigenfunction and
\begin{gather*}
\mathcal{U}_{w\left(  i\right)  }^{B}wy_{E_{k}}\zeta_{\alpha}\left(
y^{2}\right)      =2\xi_{i}\left(  \alpha\right)  wy_{E_{k}}\zeta_{\alpha
}\left(  y^{2}\right)  ,\qquad 1\leq i\leq k,\\
\mathcal{U}_{w\left(  i\right)  }^{B}wy_{E_{k}}\zeta_{\alpha}\left(
y^{2}\right)      =\left(  2\xi_{i}\left(  \alpha\right)  -1\right)
wy_{E_{k}}\zeta_{\alpha}\left(  y^{2}\right)  ,\qquad k<i\leq3.
\end{gather*}
Def\/ine $\beta$ as before ($\beta_{i}=2\alpha_{i}+1$ for $1\leq i\leq k$ and
$\beta_{i}=2\alpha_{i}$ for $k<i\leq3$) then the label for the polynomial
$wy_{E_{k}}\zeta_{\alpha}\left(  y^{2}\right)  $ is $w\beta$ (recall $\left(
w\beta\right)  _{i}=\beta_{w^{-1}\left(  i\right)  }$). Denote%
\[
p_{w\beta}\left(  y\right)  :=wy_{E_{k}}\zeta_{\alpha}\left(  y^{2}\right)  .
\]
This def\/ines a polynomial $p_{\gamma}$ for any $\gamma\in
\mathbb{N}
_{0}^{3}$. The norm of $wy_{E_{k}}\zeta_{\alpha}\left(  y^{2}\right)  $ is the
same as that of $y_{E_{k}}\zeta_{\alpha}\left(  y^{2}\right)  $ since any
$w\in\mathcal{S}_{3}$ acts as an isometry for $\left\langle \cdot
,\cdot\right\rangle _{\kappa}$. Suppose $E,E^{\prime}\subset\left\{
1,2,3\right\}  $ and $E\neq E^{\prime}$ and $f,g\in\mathcal{P}^{\left(
3\right)  }$ then $\left\langle y_{E}f\left(  y^{2}\right)  ,y_{E^{\prime}
}g\left(  y^{2}\right)  \right\rangle _{\kappa}=0$. The root system $D_{3}$ is
an orbit of the subgroup of diagonal elements of $B_{3}$ (isomorphic to $\mathbb{Z}
_{2}^{3}$). Denote the sign change $y_{i}\longmapsto-y_{i}$ by $\sigma_{i}$
for $1\leq i\leq3$. From the $B_{3}$ results we have $\sigma_{i}%
\mathcal{D}_{j}^{B}=\mathcal{D}_{j}^{B}\sigma_{i}$ for $1\leq i,j\leq3$ and
this implies $\left\langle y_{E}f\left(  y^{2}\right)  ,y_{E^{\prime}}g\left(
y^{2}\right)  \right\rangle _{\kappa}=\left\langle \sigma_{i}y_{E}f\left(
y^{2}\right)  ,\sigma_{i}y_{E^{\prime}}g\left(  y^{2}\right)  \right\rangle
_{\kappa}=-\left\langle y_{E}f\left(  y^{2}\right)  ,y_{E^{\prime}}g\left(
y^{2}\right)  \right\rangle _{\kappa}$for any $i\in\left(  E\backslash
E^{\prime}\right)  \cup\left(  E^{\prime}\backslash E\right)  $ (the symmetric
dif\/ference). Thus $\left\{  p_{\gamma}:\gamma\in
\mathbb{N}
_{0}^{3}\right\}  $ is an orthogonal basis for~$\left\langle \cdot
,\cdot\right\rangle _{\kappa}$.

We consider the $\mathcal{S}_{4}$-invariant polynomials: they are generated by
$y_{0}$, $\sum_{i=1}^{3}y_{i}^{2}$, $y_{1}y_{2}y_{3}$, $\sum_{i=1}
^{3}y_{i}^{4}.\!$ Any invariant is a sum of terms of the form $y_{0}^{n}\left(
y_{1}y_{2}y_{3}\right)  ^{s}f\left(  y^{2}\right)  $ where $n\in
\mathbb{N}
_{0},s=0$ or $1$, and $f$ is a symmetric polynomial in three variables. For
now consider only polynomials in $\left\{  y_{1},y_{2},y_{3}\right\}  $. Let
$\lambda\in \mathbb{N}
_{0}^{3,+}$, then there are two corresponding simultaneous eigenfunctions of
$\sum_{i=1}^{3}\left(  \mathcal{U}_{i}^{B}\right)  ^{n}$ (it suf\/f\/ices to take
$n=1,2,3$ to generate the commutative algebra of $\mathcal{S}_{4}$-invariant
operators). From \cite[Theorem~8.5.10]{DX} let%
\begin{gather}
A_{\lambda}     =\#\left\{  \alpha:\alpha^{+}=\lambda\right\}  \frac{\left(
3\kappa+1\right)  _{\lambda}h\left(  \lambda,1\right)  }{\mathcal{E}%
_{1}\left(  \lambda^{R}\right)  h\left(  \lambda,\kappa+1\right)
},\label{jnorm}\\
F_{\lambda}^{0}\left(  x\right)      =j_{\lambda}\left(  y^{2}\right)
,\nonumber\\
\left\langle F_{\lambda}^{0},F_{\lambda}^{0}\right\rangle _{\kappa}
=2^{2\left\vert \lambda\right\vert }\left(  2\kappa+\frac{1}{2}\right)
_{\lambda}A_{\lambda},\nonumber\\
F_{\lambda}^{1}\left(  x\right)      =y_{1}y_{2}y_{3}j_{\lambda}\left(
y^{2}\right)  ,\nonumber\\
\left\langle F_{\lambda}^{1},F_{\lambda}^{1}\right\rangle _{\kappa}
=2^{2\left\vert \lambda\right\vert }\left(  2\kappa+\frac{1}{2}\right)
_{\left(  \lambda_{1}+1,\lambda_{2}+1,\lambda_{3}+1\right)  }A_{\lambda
}.\nonumber
\end{gather}
The polynomials $\left\{  F_{\lambda}^{0},F_{\lambda}^{1}:\lambda\in \mathbb{N}
_{0}^{3,+}\right\}  $ are pairwise orthogonal.

Up to now we have mostly ignored the fourth dimension, namely, the coordinate
$y_{0}$. The ref\/lection $\sigma_{0}$ along $v_{0}$ (given by $x\sigma
_{0}=x-\left(  \sum_{i=1}^{4}x_{i}\right)  v_{0}$) commutes with the
$\mathcal{S}_{4}$-action. We introduce another parameter $\kappa^{\prime}$ and
let%
\begin{gather*}
\mathcal{D}_{0}f\left(  x\right)      =\frac{1}{2}\sum_{i=1}^{4}\frac
{\partial}{\partial x_{i}}f\left(  x\right)  +\frac{\kappa^{\prime}%
}{\left\langle x,v_{0}\right\rangle }\left(  f\left(  x\right)  -f\left(
x\sigma_{0}\right)  \right)  ,\\
\mathcal{D}_{i}^{\prime}f\left(  x\right)      =\mathcal{D}_{i}f\left(
x\right)  +\frac{\kappa^{\prime}}{2\left\langle x,v_{0}\right\rangle }\left(
f\left(  x\right)  -f\left(  x\sigma_{0}\right)  \right)  .
\end{gather*}
The operators $\left\{  \mathcal{D}_{i}^{\prime}:1\leq i\leq4\right\}  $ are
the Dunkl operators for the group $W=\mathcal{S}_{4}\times
\mathbb{Z}
_{2}$ (the ref\/lection group generated by $\left\{  \left(  1,2\right)
,\left(  2,3\right)  ,\left(  3,4\right)  ,\sigma_{0}\right\}  $). Then
$\mathcal{D}_{0}y_{0}^{2n}=2ny_{0}^{2n-1}$ and $\mathcal{D}_{0}y_{0}%
^{2n+1}=\left(  2n+1+2\kappa^{\prime}\right)  y_{0}^{2n}$. We def\/ine the
extended pairing for polynomials
\[
\left\langle f\left(  x\right)
,g\left(  x\right)  \right\rangle _{\kappa,\kappa^{\prime}}= f\left(
\mathcal{D}_{1}^{\prime},\ldots,\mathcal{D}_{4}^{\prime}\right)  g\left(
x\right)  |_{x=0};
\]
 in terms of $y$%
\begin{gather*}
   \left\langle f_{0}\left(  y_{0}\right)  f_{1}\left(  y_{1},y_{2}%
,y_{3}\right)  ,g_{0}\left(  y_{0}\right)  g_{1}\left(  y_{1},y_{2}%
,y_{3}\right)  \right\rangle _{\kappa,\kappa^{\prime}}\\
\qquad{}  =f_{0}\left(  \mathcal{D}_{0}\right)  g_{0}\left(  y_{0}\right)
|_{y_{0}=0}\times f_{1}\left(  \mathcal{D}_{1}^{B},\ldots\right)  g_{1}\left(
y_{1},y_{2},y_{3}\right)  |_{y=0}\\
\qquad{} =f_{0}\left(  \mathcal{D}_{0}\right)  g_{0}\left(  y_{0}\right)
|_{y_{0}=0}\times\left\langle f_{1},g_{1}\right\rangle _{\kappa}.
\end{gather*}
It is easily shown by induction that for $n\in
\mathbb{N}
_{0}$
\begin{gather*}
\left\langle y_{0}^{2n},y_{0}^{2n}\right\rangle _{\kappa,\kappa^{\prime}}
=2^{2n}n!\left(  \kappa^{\prime}+\frac{1}{2}\right)  _{n},\\
\left\langle y_{0}^{2n+1},y_{0}^{2n+1}\right\rangle _{\kappa,\kappa^{\prime}}
  =2^{2n+1}n!\left(  \kappa^{\prime}+\frac{1}{2}\right)  _{n+1}.
\end{gather*}
The direct product structure implies that $\left\{  p_{\left(  \gamma
_{1},\gamma_{2},\gamma_{3}\right)  }\left(  y\right)  y_{0}^{\gamma_{4}%
}:\gamma\in
\mathbb{N}
_{0}^{4}\right\}  $ is an orthogonal basis for~$\left\langle \cdot
,\cdot\right\rangle _{\kappa,\kappa^{\prime}}$.

\section{Hermite polynomials}

The pairing $\left\langle \cdot,\cdot\right\rangle _{\kappa,\kappa^{\prime}}$
is related to a measure on $\mathcal{\mathbb{R}}^{4}$: let $\kappa,\kappa^{\prime}\geq0$ and%
\begin{gather*}
dm\left(  x\right)      :=\left(  2\pi\right)  ^{-2}\exp\left(  -\frac{1}%
{2}\ \left\vert x\right\vert ^{2}\right)  dx,\qquad x\in\mathcal{\mathbb{R}
}^{4},\\
h\left(  x\right)      :=\prod_{1\leq i<j\leq4}\left\vert x_{i}-x_{j}%
\right\vert ^{\kappa}\left\vert y_{0}\right\vert ^{\kappa^{\prime}},\\
c_{\kappa,\kappa^{\prime}}^{-1}     :=\int_{\mathbb{R}
^{4}}h\left(  x\right)  ^{2}dm\left(  x\right)  ,\\
d\mu_{\kappa,\kappa^{\prime}}\left(  x\right)     :=c_{\kappa,\kappa^{\prime
}}h\left(  x\right)  ^{2}dm\left(  x\right)  .
\end{gather*}
In fact
\[
c_{\kappa,\kappa^{\prime}}^{-1}=2^{\kappa^{\prime}}\frac{\Gamma\left(
\kappa^{\prime}+\frac{1}{2}\right)  \Gamma\left(  2\kappa+1\right)
\Gamma\left(  3\kappa+1\right)  \Gamma\left(  4\kappa+1\right)  }
{\Gamma\left(  \frac{1}{2}\right)  \Gamma\left(  \kappa+1\right)  ^{3}}.
\]
The integral is a special case of the general formula (any suitably integrable
function $f$ on $\mathbb{R}$):
\begin{gather*}
 \left(  2\pi\right)  ^{-N/2}\int_{\mathbb{R}
^{N}}\prod_{1\leq i<j\leq N}\left\vert x_{i}-x_{j}\right\vert ^{2\kappa
}f\left(  \sum_{i=1}^{N}x_{i}\right)  \exp\left(  -\frac{1}{2}\ \left\vert
x\right\vert ^{2}\right)  dx\\
\qquad{}  =\frac{1}{\sqrt{2\pi}}\int_{-\infty}^{\infty}f\left(  t\sqrt{N}\right)
e^{-t^{2}/2}dt\cdot\prod_{j=2}^{N}\frac{\Gamma\left(  j\kappa+1\right)
}{\Gamma\left(  \kappa+1\right)  };
\end{gather*}
this follows from the Macdonald--Mehta--Selberg integral for $\mathcal{S}_{N}$
and the use of an orthogonal coordinate system for $\mathbb{R}
^{N}$ in which $\sum_{i=1}^{N}x_{i}/\sqrt{N}$ is one of the coordinates. The
Laplacian is $\Delta_{h}:=\sum_{i=1}^{4}\left(  \mathcal{D}_{i}^{\prime
}\right)  ^{2}=\sum_{i=1}^{3}\left(  \mathcal{D}_{i}^{B}\right)
^{2}+\mathcal{D}_{0}^{2}$. Also set $\Delta_{B}:=\sum_{i=1}^{3}\left(
\mathcal{D}_{i}^{B}\right)  ^{2}$. Then for $f,g\in\mathcal{P}$ \cite[Theorem~5.2.7]{DX}
\[
\left\langle f,g\right\rangle _{\kappa,\kappa^{\prime}}=\int_{\mathbb{R}
^{4}}\left(  e^{-\Delta_{h}/2}f\left(  x\right)  \right)  \left(
e^{-\Delta_{h}/2}g\left(  x\right)  \right)  d\mu_{\kappa,\kappa^{\prime}
}\left(  x\right)  .
\]
The orthogonal basis elements $p_{\gamma}\left(  y\right)  y_{0}^{n}$
($\gamma\in\mathbb{N}
_{0}^{3},n\in
\mathbb{N}
_{0}$) are transformed to orthogonal polynomials in $L^{2}\left(
\mathbb{R}
^{4},\mu_{\kappa,\kappa^{\prime}}\right)  $ under the action of $e^{-\Delta
_{h}/2}=\sum_{n=0}^{\infty}\frac{1}{n!}\left(  -\frac{1}{2}\right)  ^{n}%
\Delta_{h}^{n}$ (only f\/initely many terms are nonzero when acting on a
polynomial). We have%
\[
e^{-\Delta_{h}/2}\left(  p_{\gamma}\left(  y\right)  y_{0}^{n}\right)
=\left(  e^{-\Delta_{B}/2}p_{\gamma}\left(  y\right)  \right)  \left(
e^{-\mathcal{D}_{0}^{2}/2}y_{0}^{n}\right)  .
\]
Then for $n\in
\mathbb{N}
_{0}$
\begin{gather*}
e^{-\mathcal{D}_{0}^{2}/2}y_{0}^{2n}     =\left(  -2\right)  ^{n}%
n!L_{n}^{\kappa^{\prime}-\frac{1}{2}}\left(  \frac{y_{0}^{2}}{2}\right)  ,\\
e^{-\mathcal{D}_{0}^{2}/2}y_{0}^{2n+1}     =\left(  -2\right)  ^{n}%
n!y_{0}L_{n}^{\kappa^{\prime}+\frac{1}{2}}\left(  \frac{y_{0}^{2}}{2}\right)
.
\end{gather*}
Recall the Laguerre polynomials $\left\{  L_{n}^{a}\left(  t\right)  :n\in
\mathbb{N}
_{0}\right\}  $ are the orthogonal polynomials for the measure $t^{a}e^{-t}dt$
on $\left\{  t:t\geq0\right\}  $ with $a>-1$, and%
\[
L_{n}^{a}\left(  t\right)  =\frac{\left(  a+1\right)  _{n}}{n!}\sum_{i=0}%
^{n}\frac{\left(  -n\right)  _{i}}{\left(  a+1\right)  _{i}}\frac{t^{i}}{i!}.
\]
The result of applying $e^{-\Delta_{B}/2}$ to a polynomial $x_{E_{k}}%
\zeta_{\alpha}\left(  y^{2}\right)  $ is a complicated expression involving
some generalized binomial coef\/f\/icients (see \cite[Proposition~9.4.5]{DX}). For the
symmetric cases $j_{\lambda}\left(  y^{2}\right)  $ and $y_{1}y_{2}%
y_{3}j_{\lambda}\left(  y^{2}\right)  ,\lambda\in
\mathbb{N}
_{0}^{3,+}$ these coef\/f\/icients were investigated by Lassalle \cite{L} and
Okoun\-kov and Olshanski \cite[equation~(3.2)]{OO}; in the latter paper there is an
explicit formula.

Finally we can use our orthogonal basis to analyze a modif\/ication of the
type-$A$ quantum Calogero--Sutherland model with four particles on a line and
harmonic conf\/inement. By resca\-ling, the Hamiltonian (with exchange terms) can
be written as:
\[
\mathcal{H}=-\Delta+\frac{\left\vert x\right\vert ^{2}}{4}+2\kappa\sum_{1\leq
i<j\leq4}\frac{\kappa-\left(  i,j\right)  }{\left(  x_{i}-x_{j}\right)  ^{2}%
}+\frac{4\kappa^{\prime}\left(  \kappa^{\prime}-\sigma_{0}\right)  }{\left(
x_{1}+x_{2}+x_{3}+x_{4}\right)  ^{2}}.
\]
When this is applied to a $W$-invariant the ref\/lections $\left(  i,j\right)  $
and $\sigma_{0}$ are replaced by the scalar~$1$. We combine the type-$B$
results from \cite[Section~9.6.5]{DX} (setting $\kappa^{\prime}=0$ in the
formulae) with simple $\mathbb{Z}
_{2}$ calculations. The nonnormalized base state is
\[
\psi_{0}\left(  x\right)  :=\prod_{1\leq i<j\leq4}\left\vert x_{i}
-x_{j}\right\vert ^{\kappa}\left\vert y_{0}\right\vert ^{\kappa^{\prime}}
\exp\left(  -\frac{1}{4}\left\vert x\right\vert ^{2}\right)  .
\]
Then
\[
\psi_{0}^{-1}\mathcal{H\psi}_{0}=-\Delta_{B}-\mathcal{D}_{0}^{2}+\sum
_{i=0}^{3}y_{i}\frac{\partial}{\partial y_{i}}+6\kappa+\kappa^{\prime}+2.
\]
This operator has polynomial eigenfunctions and the eigenvalues are the energy
levels of the associated states. From \cite[Section~9.6.5]{DX} we have%
\[
e^{-\Delta_{B}/2}\sum_{i=1}^{3}\mathcal{U}_{i}^{B}e^{\Delta_{B}/2}=-\Delta
_{B}+\sum_{i=1}^{3}y_{i}\frac{\partial}{\partial y_{i}}+6\kappa+3,
\]
and by direct calculations%
\begin{gather*}
\mathcal{D}_{0}^{2}=\frac{\partial^{2}}{\partial y_{0}^{2}}+\frac
{2\kappa^{\prime}}{y_{0}}\frac{\partial}{\partial y_{0}}-\kappa^{\prime}
\frac{1-\sigma_{0}}{y_{0}^{2}},\\
e^{-\mathcal{D}_{0}^{2}/2}\left(  \mathcal{D}_{0}y_{0}-\kappa^{\prime}
\sigma_{0}\right)  e^{\mathcal{D}_{0}^{2}/2}=-\mathcal{D}_{0}^{2}+y_{0}
\frac{\partial}{\partial y_{0}}+\kappa^{\prime}+1.
\end{gather*}
Combine these results:
\[
\psi_{0}^{-1}\mathcal{H\psi}_{0}=e^{-\Delta_{h}/2}\left(  \sum_{i=1}%
^{3}\mathcal{U}_{i}^{B}+\mathcal{D}_{0}y_{0}-\kappa^{\prime}\sigma
_{0}-2\right)  e^{\Delta_{h}/2}.
\]
Thus $\left(  e^{-\Delta_{h}/2}\left(  p_{\gamma}\left(  y\right)  y_{0}%
^{n}\right)  \right)  \psi_{0}$ is an eigenfunction of $\mathcal{H}$ for each
$\gamma\in
\mathbb{N}
 _{0}^{3}$, $n\in \mathbb{N}
_{0}$. It suf\/f\/ices to consider $y_{E_{k}}\zeta_{\alpha}\left(  y^{2}\right)
y_{0}^{n}$. We have%
\begin{gather*}
\left(  \mathcal{D}_{0}y_{0}-\kappa^{\prime}\sigma_{0}\right)  y_{0}^{2n}
=\left(  \left(  2n+1+2\kappa^{\prime}\right)  -\kappa^{\prime}\right)
y_{0}^{2n},\\
\left(  \mathcal{D}_{0}y_{0}-\kappa^{\prime}\sigma_{0}\right)  y_{0}^{2n+1}
   =\left(  \left(  2n+2\right)  +\kappa^{\prime}\right)  y_{0}^{2n},\\
\left(  \mathcal{D}_{0}y_{0}-\kappa^{\prime}\sigma_{0}\right)  y_{0}^{n}
=\left(  n+1+\kappa^{\prime}\right)  y_{0}^{n}.
\end{gather*}
Furthermore $\sum_{i=1}^{3}\mathcal{U}_{i}^{B}\left(  y_{E_{k}}\zeta_{\alpha
}\left(  y^{2}\right)  \right)  =\left(  2\sum_{i=1}^{3}\xi_{i}\left(
\alpha\right)  -\left(  3-k\right)  \right)  y_{E_{k}}\zeta_{\alpha}\left(
y^{2}\right)  $; the eigenvalue is $\left(  2\left\vert \alpha\right\vert
+k\right)  +6\kappa+3=\left\vert \beta\right\vert +6\kappa+3$ (where
$\beta_{i}=2\alpha_{i}+1$ for $1\leq i\leq k$ and $\beta_{i}=2\alpha_{i}$ for
$k<i\leq3$). The energy level for $\left(  e^{-\Delta_{h}/2}\left(  p_{\beta
}\left(  y\right)  y_{0}^{n}\right)  \right)  \psi_{0}$ is $\left\vert
\beta\right\vert +n+6\kappa+\kappa^{\prime}+2$. Observe the degeneracy of the
energy levels; only the total degree $\left\vert \beta\right\vert +n$ appears.
The (nonnormalized) $W$-invariant eigenfunctions are ($\lambda\in%
%TCIMACRO{\U{2115} }%
%BeginExpansion
\mathbb{N}
%EndExpansion
_{0}^{3}$)%
\begin{align*}
&  \left(  e^{-\Delta_{B}/2}\left(  j_{\lambda}\left(  y^{2}\right)  \right)
L_{n}^{\kappa^{\prime}-1/2}\left(  \frac{y_{0}^{2}}{2}\right)  \right)
\psi_{0}\left(  x\right)  ,\\
&  \left(  e^{-\Delta_{B}/2}\left(  y_{1}y_{2}y_{3}j_{\lambda}\left(
y^{2}\right)  \right)  L_{n}^{\kappa^{\prime}-1/2}\left(  \frac{y_{0}^{2}}%
{2}\right)  \right)  \psi_{0}\left(  x\right)  .
\end{align*}
The $L^{2}$-norms can be found by using equation (\ref{jnorm}).

In conclusion, we have found an unusual basis for polynomials which allowed an
extra parameter in the action of $\mathcal{S}_{4}$ on $%
%TCIMACRO{\U{211d} }%
%BeginExpansion
\mathbb{R}
%EndExpansion
^{4}$. This exploited the fact that $v_{0}^{\bot}$ has an orthogonal basis
which together with its antipodes forms an $\mathcal{S}_{4}$-orbit. The
pairing $\left\langle \cdot,\cdot\right\rangle _{\kappa}$ has an analog for
each ref\/lection group and weight function. We are left with the interesting
problem of how to construct orthogonal bases for groups not of type $A$ or $B$.

\pdfbookmark[1]{References}{ref}
\LastPageEnding

\end{document}